\documentclass[a4paper, 12pt]{article}
\usepackage{a4,amsmath,amsthm,amssymb,amsfonts}
\usepackage{pdflscape}
\usepackage{verbatim}
\newcommand{\R}{{\mathbb R}}

\newcommand{\Z}{{\mathbb Z}}
\newcommand{\N}{{\mathbb N}}

\newcommand{\be}{\begin{equation}}
\newcommand{\ee}{\end{equation}}

\newcommand \I   {\mathrm{i}}

\newcommand \ra {\rightarrow}

\newcommand   \al  \alpha
\newcommand   \bt  \beta
\newcommand   \gm  \gamma
\newcommand   \dl  \delta
\newcommand   \lm  \lambda
\newcommand   \om  \omega
\newcommand   \Om  \Omega
\newcommand   \sg  \sigma
\newcommand   \eps \varepsilon
\renewcommand \phi \varphi

\newcommand \induces       \leadsto
\newcommand \contradiction \lightning
\newcommand \ctrdct        \contradiction

\newcommand \ds            \displaystyle
\newcommand \lbar          \overline

\newtheoremstyle{thm}
  {\baselineskip}
  {\baselineskip}
  {\itshape}
  {}
  {\bfseries}
  {.}
  {.5em}
  {}

\newtheoremstyle{others}
  {\baselineskip}
  {\baselineskip}
  {\upshape}
  {}
  {\bfseries}
  {.}
  {.5em}
  {}

\theoremstyle{thm}

\newtheorem* {satz*}                {Satz}

\newtheorem* {corollary*}           {Corollary}

\theoremstyle{others}

\newtheorem* {remark*}              {remark}

\newtheorem* {definition*}          {Definition}

\newtheorem* {example*}             {Example}

\newtheorem* {notation*}            {Bezeichnung}

\numberwithin{theorem}{section}
\numberwithin{equation}{section}

\begin{document}
\title{A sharp version of Ehrenfest's theorem \\ for general self-adjoint operators}

\author{Gero Friesecke and Bernd Schmidt \\ Center for Mathematics, TU Munich \\
        {\tt gf@ma.tum.de}, {\tt schmidt@ma.tum.de}}

\date{December 2, 2009}

\maketitle

\begin{abstract} We prove the Ehrenfest theorem of quantum mechanics under
sharp assumptions on the operators involved. 
\begin{center}
$\mbox{\tiny AMS Subject classification: 81Qxx, 47B25, 47D08}$ \\
\end{center}
\end{abstract}

{\bf Introduction.} $\;$ While we know a great deal about analytic properties
of the \linebreak commutative algebra generated
by a single unbounded self-adjoint operator 
(via \linebreak a continuous or measurable functional calculus), 
much less is known mathematically
about the non-commutative algebra generated by more than one such operator. 

This note is concerned with Ehrenfest's theorem in quantum mechanics, which naturally
involves pairs of noncommuting unbounded self-adjoint
operators. For a quantum system which evolves via the
time-dependent Schr\"odinger equation 
\be \label{SE}
          i\frac{d}{dt}\psi(t)=H\psi(t),
\ee
Ehrenfest's theorem asserts that
the time evolution of the quadratic forms
$$
       \langle\psi(t),\, A\psi(t)\rangle = 
       \langle\psi(0), \, e^{itH}Ae^{-itH}\psi(0)\rangle
$$
is governed by the equation 
\be \label{eh3}
      \frac{d}{dt}\langle A\rangle_{\psi(t)}=i\langle [H,A]\rangle_{\psi(t)}.
\end{equation}
Here $\psi(t)$ belongs to a complex Hilbert space $X$, $H$ and $A$
are unbounded self-adjoint operators on $X$, $[H,A]$ denotes the commutator
$HA-AH$, and $\langle A \rangle_{\phi}$ stands for the quadratic form $\langle\phi, A\phi\rangle$. 
Of course, at the moment the terms in the above equation are only formal, since it is not clear whether and where
the required operator compositions are well defined. 

Physically, 
$\psi(t)$ describes the state of the quantum system at time $t$, 
$H$ is the Hamiltonian of the system (obtained, e.g., by
quantizing an underlying classical Hamiltonian), $A$ is an observable, and the value of its quadratic form gives the expected
value of the observable in a measurement.
In Ehrenfest's original work \cite{Ehrenfest}, $H$ was a Schr\"odinger operator
of form $-\Delta+V(x)$ and $A$ a component of the position or momentum operator. In this case 
the left hand side of \eqref{eh3} describes the expected value of a position or momentum measurement, 
and the right hand side of \eqref{eh3} turns out to be the expected value of velocity
respectively force, linking quantum mechanics to classical mechanics.
Rigorous versions are not difficult to obtain in the context of classical, rapidly decaying solutions for smooth potentials \cite{BlankExnerHavlicek} or when $A$ is relatively bounded with respect to $H$, 
and play an important role in the analysis of space-time behaviour of Schr\"odinger wavefunctions
\cite{Hunziker, RS} and scattering theory \cite{Sigal, Graf,
Derezinski, DG}. For a rigorous treatment of non-relatively-bounded $A$ and applications to 
atomic and molecular systems with Coulomb interactions see \cite{FK}. 

Our goal in this note is to prove the following sharp version of Ehrenfest's theorem for general
self-adjoint operators. 
\\[4mm]
{\bf Theorem 1} {\it
\label{T:1} (Sharp Ehrenfest theorem) Let $A \, : \, D(A)\to X$, $H\, : \, D(H)\to X$ be densely defined
self-adjoint operators on a Hilbert space $X$. If $e^{-itH}$ leaves $D(A)\cap D(H)$ invariant for all $t$, then for any $\psi_0\in D(A)\cap D(H)$, 
the expected value $\langle A \rangle_{\psi(t)}$, $\psi(t):=e^{-itH}\psi_0$, is 
continuously differentiable with respect to $t$ and satisfies eq. (\ref{eh3}), the right hand side 
being understood in the inner product sense
\be \label{formdef}
    \langle [H,A]\rangle_{\psi}:= \langle H\psi,A\psi\rangle-\langle A\psi,H\psi\rangle \;\; 
    \mbox{($\psi\in D(A)\cap D(H)$)}.
\ee
} 
\vspace*{-2mm}

\noindent
(Recall that a linear operator $A : D(A) \to X$ defined on a subspace $D(A)$ of a Hilbert space $X$ is called hermitean when $\langle A \psi, \varphi \rangle = \langle \psi, A \varphi \rangle$ for all $\psi, \varphi \in D(A)$, and self-adjoint if in addition the domain $D(A)$ is dense and in a suitable sense maximal, see e.g.\ \cite{RSI}. Self-adjoint operators $H$ generate a one-parameter unitary group $e^{-itH}$ \cite{RSI}.) 
Note that the hypothesis above that $\psi(t)=e^{-itH}\psi_0$ belongs to $D(A)$ for all $t$ is necessary for the element
$A\psi(t)$ appearing in the inner product 
in (\ref{eh3}) to be well defined. Thus Ehrenfest's theorem always holds when it makes sense in an operator-theoretic setting. 
(Conceivable abstract extensions of the expressions in (\ref{eh3}) to quadratic form domains are not investigated here.)

>From an analysis
point of view, an interesting aspect of the result is the automatic higher regularity of the expected value $\langle A \rangle_{\psi(t)}$ in
time; a priori it is not even clear that $t\mapsto\langle A \rangle_{\psi(t)}$ is continuous.

Self-adjointness of $A$ cannot be weakened to hermiteanity; otherwise,
even continuity of $t\mapsto\langle A \rangle_{\psi(t)}$ can fail. See the counterexample at the end of this note.
This phenomenon is not obvious, since only $H$ is required to generate a unitary
group in order for the expected values appearing in (\ref{eh3}) to be well defined. 
\\[6mm]
{\bf Proof of Theorem 1.}
Starting point is the following result of \cite{FK}, which involves a weaker assumption on
the operator $A$ (densely defined hermitean instead of self-adjoint) but a stronger assumption
on the ``coupling'' between the unitary group $e^{-itH}$ and the operator $A$ (see (H3)
below):
\\[4mm]
{\bf Proposition 1} \cite{FK} {\it
\label{T:2} Let $H$ and $A$ be two densely defined linear operators on $X$
such that: 
\\[1mm]
{\rm (H1)} $H:D(H)\ra X$ is selfadjoint, $A:D(A)\ra X$ is hermitean \\
{\rm (H2)} ${\rm e}^{-\I tH}$ leaves $D(A)\cap D(H)$ invariant for all $t\in\R$ \\
{\rm (H3)} For any $\psi_0\in D(A)\cap D(H)$, $\sup_{t\in I}||Ae^{-itH}\psi_0||<\infty$ for $I\subset\R$ bounded.
\\[1mm]
Then for $\psi_0\in D(A)\cap D(H)$, the expected value $\langle A \rangle_{\psi(t)}$, $\psi(t):=e^{-itH}\psi_0$, is 
continuously differentiable with respect to $t$ and satisfies eq. (\ref{eh3}), the right hand side being understood in the
sense of (\ref{formdef}).
} 
\\[4mm]
%
%
%
For completeness we sketch the (simple) proof. Consider the difference quotient 
\begin{align*} 
    &\frac{\langle \psi(t\!+\! h),\, A\psi(t\!+\! h)\rangle - \langle \psi(t), \, A\psi(t)\rangle}{h} \\
    &=\Bigl\langle A\psi(t\!+\! h), \, \frac{\psi(t\!+\! h)-\psi(t)}{h} \Bigr\rangle
    +
    \Bigl\langle \frac{\psi(t\!+\! h)-\psi(t)}{h}, \, A\psi(t)\Bigr\rangle.
\end{align*} 
As $h\to 0$, the second term goes to $i\langle H\psi(t),\, A\psi(t)\rangle$, since $\frac{\psi(t+h)-\psi(t)}{h}\to -iH\psi(t)$
strongly due to the strong differentiability of the unitary group $\{e^{-itH}\}_{t\in\R}$ on $D(H)$. The first term
goes to $-i\langle A\psi(t),\, H\psi(t)\rangle$ since the first factor $A\psi(t+h)$, being bounded by (H3), converges weakly
up to subsequences to a limit $f\in X$, which must equal $A\psi(t)$ since
\begin{align*}
  \langle f, \, \phi\rangle = \lim_{h_j\to 0}\langle A\psi(t\!+\! h_j), \,\phi\rangle = \lim_{h_j\to 0}
  \langle\psi(t\!+\! h_j), A\phi\rangle
  = \langle \psi(t), A\phi \rangle
  = \langle A\psi(t), \phi\rangle
\end{align*}
for $\phi$'s belonging to the dense subset $D(A)$. Hence the whole sequence $A\psi(t+h)$ tends weakly to $A\psi(t)$,
and so $\langle A \rangle_{\psi(t)}$ is differentiable with respect
to $t$ and satisfies eq. (\ref{eh3}). Continuity of the derivative follows by combining the above
weak convergence with the strong convergence of $H\psi(t+h)$ to $H\psi(t)$. 

Theorem {1} is an immediate consequence of Proposition 1 and the following
\\[4mm]
{\bf Proposition 2} {\it 
\label{AandHselfadj} Let $H$ and $A$ be two densely defined operators on $X$ such that \\[1mm]
{\rm (H1')} $H \, : \, D(H)\to X$ self-adjoint, $A\, : \, D(A)\to X$ hermitean and closed \\
{\rm (H2)} $T(t):=e^{-itH}$ leaves $D(A)\cap D(H)$ invariant for all $t\in\R$. \\[1mm]
Then for any $\psi_0\in D(A)\cap D(H)$, $\sup_{t\in I}||A T(t)\psi_0||<\infty$ for $I\subset\R$ bounded.
} 
\\[4mm]
{\bf Proof of Proposition 2.} Throughout the proof we fix $\psi_0\in D(A)\cap D(H)$. 
\\[1mm]
{\bf 1.} We claim that there exists a nonempty interval $I=(a,b)$ such that $\sup_{t\in I}$ $||A T(t)\psi_0||<\infty$.

This can be established with the help of the Baire category theorem, as follows. Consider the sets
$E_n:=\{t\in\R\, : \, ||AT(t)\psi_0||\le n\}$, $n\in\N$. Because $T(t)\psi_0\in D(A)$ for all $t$, we have 
$\cup_{n\in\N}E_n=\R$. Moreover the $E_n$ are closed. Namely, suppose $t_j\in E_n$, $t_j\to t$. Then 
$||AT(t_j)\psi_0||\le n$, so $\{AT(t_j)\psi_0\}$ is a bounded sequence in $X$. Hence for a subsequence, 
again denoted $t_j$, $AT(t_j)\psi_0\rightharpoonup\chi$ weakly in $X$, and by weak lower semi-continuity of
the norm, $||\chi||\le \liminf_{j\to\infty}||AT(t_j)\psi_0||$, whence $||\chi||\le n$. To infer that $E_n$
is closed, it suffices to show that $\chi=AT(t)\psi_0$. But for all $\eta\in D(A)$, using the weak convergence
of $AT(t_j)\psi_0$, the hermiteanity of $A$, and the continuity of $T(t)\psi_0$ in $t$, 
$$
  \langle\eta,\chi\rangle = \lim_{j\to\infty}\langle\eta, AT(t_j)\psi_0\rangle 
 =\lim_{j\to\infty}\langle A\eta, T(t_j)\psi_0\rangle = \langle A\eta, T(t)\psi_0\rangle 
 =\langle\eta, AT(t)\psi_0\rangle.
$$
Since $D(A)$ is dense, $AT(t)\psi_0=\chi$, as asserted. By the Baire category theorem, the countable union 
of the $E_n$ can't be nowhere dense; so the closure of some $E_n$ must contain an open set. Since this $E_n$
is closed, it must contain an open set, and the assertion of Step 1 follows. We remark that this step did not use
the closedness of the operator $A$.
\\[2mm]
{\bf 2.} We claim that there exists an open interval containing $0$ such that $\sup_{t\in I}$ $||A T(t)\psi_0||<\infty$.

To see this, the idea is to use the group property of $T(t)$ as well as an appropriate norm on $D(A)\cap D(H)$ with
respect to which $AT(t)$ is a continuous map into $X$. Let $t_0:=\frac{a+b}{2}$ be the midpoint of the
interval from Step 1, and let $t\in(a,b)$. Then
\be \label{generalbound}
  ||AT(-t_0+t)\psi_0|| = || AT(-t_0)T(t)\psi_0|| \le \sup_{\eta\in D(A)\cap D(H)\backslash\{0\}}
         \frac{||AT(-t_0)\eta||}{||\eta||_E} \; \cdot \; ||T(t)\psi_0||_E,
\ee
for any norm $||\cdot||_{E}$ on $D(A)\cap D(H)=:E$. The art now is to choose the norm in such a way that the first
factor on the right is finite, and the second factor is uniformly bounded for $t\in (a,b)$. A natural candidate for
such a norm is $||\cdot||_E := ||\cdot|| + ||A\cdot|| + ||H\cdot||$. By Step 1 and the fact that $||T(t)\psi_0||$ and
$||HT(t)\psi_0||$ are conserved, the second factor satisfies
\begin{eqnarray} 
  \sup_{t\in(a,b)}||T(t)\psi_0||_{E} & = & \sup_{t\in (a,b)} \Bigl(||T(t)\psi_0|| + ||AT(t)\psi_0|| + ||HT(t)\psi_0||\Bigr) 
  \nonumber \\
  & = & ||\psi_0|| + ||H\psi_0|| + \sup_{t\in (a,b)}||AT(t)\psi_0|| < \infty. \label{star}
\end{eqnarray}
To see that the first factor on the RHS of (\ref{generalbound}) is finite is less easy. 

First, we claim that $E$
with the above norm is a Banach space. This is a variant of the familiar fact that the domain $D(A)$ of a closed
operator endowed with the graph norm $||\psi||+||A\psi||$ is a Banach space, and is straightforward to verify 
using the closedness of $A$ and $H$. 


Next, we claim that $B:=AT(-t_0)$ is a closed operator from $E$ to $X$. (Clearly $B$ is defined on all of $E$,
because $T(-t_0)$ maps $D(A)$ into $D(A)$ by (H2).) To prove this, we have to show that if 
$\psi_j\in E$, $\psi_j\to\psi$ in $E$ and $B\psi_j\to\eta$ in $X$, then $\psi\in D(B)$ and $B\psi=\eta$. 
Since $B$ is defined on all of $E$, the first condition is satisfied. To show that $B\psi=\eta$, we use that
$||\psi_j-\psi||_X\le ||\psi_j-\psi||_E$, so $\psi_j\to\psi$ in $X$. By the boundedness of $T(-t_0)\, : \, X\to X$,
it follows that $T(-t_0)\psi_j\to T(-t_0)\psi$ in $X$. Because $A T(-t_0)\psi_j\to\eta$ in $X$ and $A$ is closed,
$T(-t_0)\psi\in D(A)$ and $AT(-t_0)\psi=\eta$. Thus $B\psi=\eta$, as was to be shown. 

Finally, we appeal to the closed graph theorem: $B\, : \, E\to X$ is a closed linear operator between Banach spaces,
and must hence be bounded. Consequently
\be \label{twostar}
  ||| AT(-t_0) |||_{L(E,X)} := \sup_{\eta\in E\backslash\{0\}} \frac{ ||AT(-t_0)\eta||}{||\eta||_E } < \infty 
  \mbox{ for all }t_0\in\R.
\ee
Together with (\ref{generalbound}), (\ref{star}), this establishes the assertion of Step 2. 
\\[2mm]
{\bf 3.} We claim that on any bounded interval $I\subset\R$, $\sup_{t\in I}||A T(t)\psi_0||<\infty$. 

By Step 2, there exist $\epsilon>0$, $C>0$ such that $||AT(t)\psi_0||\le C$ for all $|t|\le \epsilon$. Next, consider
$\eps\le |t|\le 2\eps$. For each such $t$, $T(t)$ can be written as $T(\eps)T(s)$ or $T(-\eps)T(s)$ for some $|s|\le\eps$.
Clearly for $\psi_0\in D(A)\cap D(H)=E$
\be
   || AT(\pm\epsilon) T(s)\psi_0|| \le |||AT(\pm\eps)|||_{L(E,X)} \sup_{|s|\le\epsilon} ||T(s)\psi_0||_E 
                                   \le |||AT(\pm\eps)|||_{L(E,X)} \, C,
\ee
the above operator norm being finite due to (\ref{twostar}). 
Hence
$$
   ||AT(t)\psi_0|| \le \underbrace{\max\Bigl\{ |||AT(\eps)|||_{L(E,X)}, \, |||AT(-\eps)|||_{L(E,X)}\Bigr\}}_{=:M} \, C 
   \mbox{ for all }|t|\le 2\epsilon.
$$
Iterating gives $||AT(t)\psi_0||\le M^{n-1}C$ for all $|t|\le n\epsilon$, for any $n\in\N$. This completes the proof
of Proposition 2.

Proposition 2 could alternatively be proved by recourse to known but nontrivial results in the theory of $C_0$-groups. One first
deduces from restriction properties of $C_0$-groups \cite[Proposition 3.2.5]{XY} 
that $T(t)$ defines a $C_0$-group on $D(A)\cap D(H)$ with respect to the graph norm of the pair of $A$ and $H$. This, however,
requires a considerable amount of machinery, in particular a characterization of $C_0$-groups on Banach spaces with separable
dual via weak Borel measurability. One then appeals to known locally uniform bounds on norms of $C_0$-groups \cite[Proposition 3.2.2]{XY}
which are related to the fact that measurable solutions to the functional inequality $e^{t+s}\le e^{t}e^{s}$ are locally uniformly
bounded. By contrast, the proof given above only uses basic functional analysis.

Observe that in fact it suffices to assume in Theorem {1} that $A$ be a closed hermitean 
operator, neither required to be densely defined nor self-adjoint. This can be seen by 
replacing $A$ by $P A P$, where $P$ denotes the orthogonal projection from $X$ to 
$\tilde{X} := \overline{\operatorname{span}} \{ \psi(t) : t \in \R \}$, and by applying 
Proposition {2} and Proposition {1} to $\tilde{X}$ instead of ${X}$.
\\[6mm]
{\bf Counterexample.} 
The following example shows that the assumption in Ehrenfest's theorem that $A$ be
self-adjoint cannot be weakened to hermitean
(and that the assumption in Proposition 2 that $A$ be closed cannot be dropped). 

Let $X = L^2(\R)$, $H = -i \frac{d}{d x}$, so that the unitary group $e^{-i t H}$ acts by translation: 
$e^{- i t H} \psi_0(x) = \psi_0(x - t)$. The idea is that the finite linear span of a well chosen single orbit
$\{\psi(t) \, : \, t\in\R\}$, $\psi(t)=e^{-itH}\psi_0$, contains infinitely many functions $\varphi_j$ with disjoint support contained in
a single bounded interval. The $\varphi_j$ can be constructed such that $\psi(t)$ becomes ``resonant'' with 
the projectors $P_j=|\varphi_j\rangle \, \langle\varphi_j|$ at particular times $t_j$, whereas at these times it is annihilated by   
all the other projectors. One then takes $A$ to be a suitable linear combination of these projectors $P_j$ on 
$\operatorname{span} \{ \psi(t) : t \in \R \}$, extended by $0$ on the orthogonal complement. 

Let $\psi_0 \in D(H) = H^1(\R)$ be the ``tent'' function with $\psi_0(0) = \psi_0(2) = 0$, $\psi_0(1) = 1$ and such that $\psi_0$ is linear 
on each of the intervals $(-\infty, 0)$, $(0, 1)$, $(1, 2)$, $(2, \infty)$. 

{\bf 1.} If $I$ is any open subinterval of $(0, 1)$, we can choose $\varphi \in \operatorname{span} \{ \psi(t) : t \in I\} \setminus \{ 0 \}$
(meaning finite linear combinations) 
such that $\varphi \perp \operatorname{span} \{ \psi(t) : t \notin I + \Z \}$ 
in the following way. Let $t_0, t_0 + 6 \eta \in I$ and in order to achieve small support consider the second difference
$\tilde{\varphi} := \psi(t_0) - 2 \psi(t_0 + \eta) + \psi(t_0 + 2\eta)
=\psi_0(\cdot-t_0)-2\psi_0(\cdot-(t_0+\eta))+\psi_0(\cdot-(t_0+2\eta))$. 
Then $\tilde{\varphi} \ne 0$ is supported on $[t_0, t_0 + 2\eta] + \{0, 1, 2\}$. In fact, on $(0,1)$ $\tilde{\varphi}$ is the piecewise linear interpolation of $\tilde{\varphi}(t_0) = 0, \tilde{\varphi}(t_0 + \eta) = \eta, \tilde{\varphi}(t_0 + 2\eta) = 0$ and 
$\tilde{\varphi}(t + 1) = - 2 \tilde{\varphi}(t)$, $\tilde{\varphi}(t + 2) =  \tilde{\varphi}(t)$ for $t \in [0, 1]$.

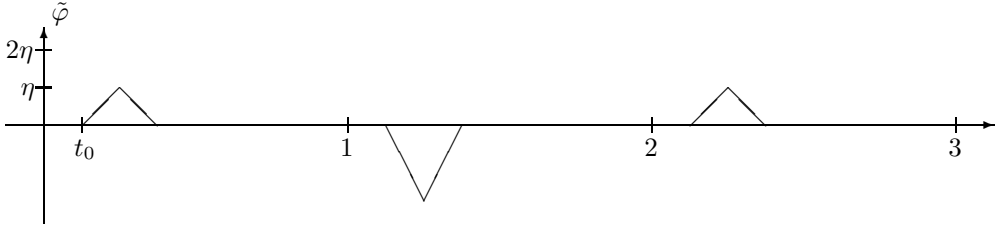
\setlength{\unitlength}{1cm}
\begin{picture}(15,3.5)
\begin{footnotesize}
\put(0,1.5){\vector(1,0){13}}
\put(0.5,0.2){\vector(0,1){2.6}}
\put(0.4,2.5){\line(1,0){.2}}
\put(0,2.4){$2\eta$}
\put(0.4,2){\line(1,0){.2}}
\put(0.2,1.9){$\eta$}
\put(0.6,2.9){$\tilde{\varphi}$}
\put(1,1.4){\line(0,1){.2}}
\put(0.9,1.1){$t_0$}

\put(4.5,1.4){\line(0,1){.2}}
\put(4.4,1.1){$1$}
\put(8.5,1.4){\line(0,1){.2}}
\put(8.4,1.1){$2$}
\put(12.5,1.4){\line(0,1){.2}}
\put(12.4,1.1){$3$}

\put(1,1.5){\line(1,1){.5}}
\put(2,1.5){\line(-1,1){.5}}
\put(5,1.5){\line(1,-2){.5}}
\put(6,1.5){\line(-1,-2){.5}}
\put(9,1.5){\line(1,1){.5}}
\put(10,1.5){\line(-1,1){.5}}

\end{footnotesize}
\end{picture}

Now to achieve orthogonality to affine functions take again a second difference, 
$\varphi := \tilde{\varphi} - 2 \tilde{\varphi}(\cdot - 2\eta) + \tilde{\varphi}(\cdot - 4\eta) = 
\psi(t_0) - 2 \psi(t_0 + \eta) - \psi(t + 2\eta) + 4 \psi(t + 3\eta) - \psi(t + 4\eta) - 2 \psi(t + 5\eta) + \psi(t + 6\eta)$. 
Then $\varphi \in \operatorname{span} \{ \psi(t) : t \in I + \Z \} \setminus \{ 0 \}$, and $\varphi$ is supported on $[t_0, t_0 + 2\eta] + \{0, 1, 2\}$ 

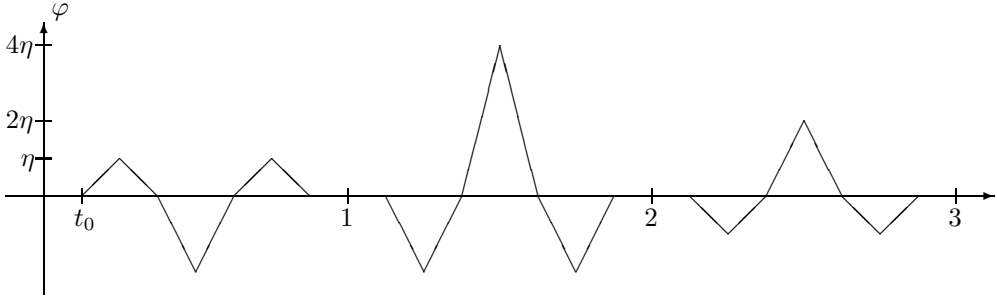
\setlength{\unitlength}{1cm}
\begin{picture}(15,4.2)
\begin{footnotesize}
\put(0,1.5){\vector(1,0){13}}
\put(0.5,0.2){\vector(0,1){3.6}}
\put(0.4,3.5){\line(1,0){.2}}
\put(0,3.4){$4\eta$}
\put(0.4,2.5){\line(1,0){.2}}
\put(0,2.4){$2\eta$}
\put(0.4,2){\line(1,0){.2}}
\put(0.2,1.9){$\eta$}
\put(0.6,3.9){$\varphi$}
\put(1,1.4){\line(0,1){.2}}
\put(0.9,1.1){$t_0$}

\put(4.5,1.4){\line(0,1){.2}}
\put(4.4,1.1){$1$}
\put(8.5,1.4){\line(0,1){.2}}
\put(8.4,1.1){$2$}
\put(12.5,1.4){\line(0,1){.2}}
\put(12.4,1.1){$3$}

\put(1,1.5){\line(1,1){.5}}
\put(2,1.5){\line(-1,1){.5}}
\put(2,1.5){\line(1,-2){.5}}
\put(3,1.5){\line(-1,-2){.5}}
\put(3,1.5){\line(1,1){.5}}
\put(4,1.5){\line(-1,1){.5}}

\put(5,1.5){\line(1,-2){.5}}
\put(6,1.5){\line(-1,-2){.5}}
\put(6,1.5){\line(1,4){.5}}
\put(7,1.5){\line(-1,4){.5}}
\put(7,1.5){\line(1,-2){.5}}
\put(8,1.5){\line(-1,-2){.5}}

\put(9,1.5){\line(1,-1){.5}}
\put(10,1.5){\line(-1,-1){.5}}
\put(10,1.5){\line(1,2){.5}}
\put(11,1.5){\line(-1,2){.5}}
\put(11,1.5){\line(1,-1){.5}}
\put(12,1.5){\line(-1,-1){.5}}


\end{footnotesize}
\end{picture}

\noindent and one easily sees that $\int_{t_0 + k}^{t_0 + 6 \eta + k} x \varphi(x) \, dx = 0$ for all $k \in \Z$. 
Now if $t \notin I + \Z$, then $\psi(t)$ is linear on each interval $(t_0 + k, t_0 + 6 \eta + k)$, so we obtain that indeed $\langle \psi(t), \varphi\rangle = 0$.

{\bf 2.} Now choose infinitely many pairwise disjoint open intervals $I_j \subset (0,1)$, $j \in \N$. Choosing vectors corresponding to the 
interval $I_j$ as described above and normalizing, we obtain an orthonormal system $(\varphi_j)$. To each $\varphi_j \in \operatorname{span} \{ \psi(t) : t \in I_j \}$ we may also choose $t_j \in I_j$ and $a_j \in \R$ such that $a_j |\langle \varphi_j, \psi(t_j) \rangle|^2 > j$. We then define the hermitean operator $A$ on $D(A) := \operatorname{span} \{ \psi(t) : t \in \R \} \oplus \operatorname{span} \{ \psi(t) : t \in \R \}^{\perp}$ by 
\begin{align*}
  A \psi := \sum_{j \in \N} a_j \langle \varphi_j, \psi \rangle \varphi_j
\end{align*}
on $\operatorname{span} \{ \psi(t) : t \in \R \}$ and zero on its orthogonal complement.
(Note that by construction the sum in this definition is in fact finite.) Then $A$ is a densely defined hermitean operator such that 
$D(A) \cap D(H)$ is invariant under $e^{-i t H}$. But $\langle A \rangle_{\psi(t)}$ is not bounded on $(0, 1)$ because on $I_j$ 
\begin{align*}
  \langle \psi(t), A \psi(t) \rangle 
  = \langle \psi(t), \sum_{j \in \N} a_j \langle \varphi_j, \psi(t) \rangle \varphi_j \rangle 
  = a_j |\langle \psi(t), \varphi_j \rangle|^2, 
\end{align*}
in particular, $\langle \psi(t_j), A \psi(t_j) \rangle > j$. 

\end{document}